\documentclass[10pt]{amsart}
\usepackage{fullpage, amsfonts,amsmath,amscd,amssymb}
\usepackage{amssymb, amsbsy, amsthm, amsmath, amstext, amsopn, verbatim,
multicol}
\usepackage{calc}
\usepackage{ifthen}

\newtheorem*{thm}{Theorem}

\newtheorem*{lem}{Lemma}

\theoremstyle{definition}

\newtheorem*{conjecture}{Conjecture}

\theoremstyle{remark}

\newtheorem*{ackn}{Acknowledgements}

\newcommand{\im}{\operatorname{im}}

\newcommand{\h}{\mathfrak{h}}

\newcommand{\Z}{\mathbb{Z}}
\newcommand{\Q}{\mathbb{Q}}

\newcommand{\C}{\mathbb{C}}

\DeclareMathOperator{\id}{id}

    \newcommand{\btheta}{\boldsymbol{\theta}}  
    
    \newcommand{\blambda}{\boldsymbol{\lambda}}
    
    \newcommand{\bmu}{\boldsymbol{\mu}}

\DeclareMathOperator{\irr}{\textsf{Irr}}

\begin{document}
\title{Calogero-Moser space, reduced rational Cherednik algebras and two-sided cells.}
\author{I.G. Gordon} \address{School of Mathematics and Maxwell Institute for Mathematical Sciences,
University of Edinburgh,
James Clerk Maxwell Building,
Kings Buildings, Mayfield Road, Edinburgh EH9 3JZ,
Scotland}
\email{igordon@ed.ac.uk}
\author{M. Martino} \address{Mathematisches Institut, Universit\"at zu K\"oln, Weyertal 86-90, D-50931 K\"oln, Germany} \email{mmartino@mi.uni-koeln.de}

\begin{abstract} We conjecture that the ``nilpotent points" of Calogero-Moser space for reflection groups are parametrised naturally by the
two-sided cells of the group with unequal parameters. The nilpotent points correspond to blocks of
restricted Cherednik algebras and we
describe these blocks in the case $G = \mu_{\ell}\wr \mathfrak{S}_n$
and show that in type $B$ our description produces an existing conjectural
description of two-sided cells.\end{abstract} \maketitle
\section{Introduction}
\subsection{} Smooth points are all alike; every singular point is singular in its own way. Calogero-Moser space associated to the symmetric group has remarkable applications in a broad range of topics; in \cite{EG}, Etingof and Ginzburg introduced a generalisation associated to any complex reflection group which has also found a variety of uses. The Calogero-Moser spaces associated to a complex reflection group, however, exhibit new behaviour: they are often singular. The nature of these singularities remains a mystery, but  their existence has been used to solve the problem of the existence of crepant resolutions of symplectic quotient singularities. The generalised Calogero-Moser spaces are moduli spaces of representations of rational Cherednik algebras and so their geometry reflects the representation theory of these algebras: smooth points correspond to irreducible representations of maximal dimension; singular points to smaller, more interesting representations. In this note we conjecture a strong link between the representations corresponding to some particularly interesting ``nilpotent points" of Calogero-Moser space and Kazhdan-Lusztig cell theory for Hecke algebras with unequal parameters. To justify the conjecture we give a combinatorial parametrisation of these points, thus answering a question of \cite{G03}, and then relate this parametrisation to the conjectures of \cite{BGIL06} on the cell theory for Weyl groups of type $B$.

\subsection{} Let $W$ be a complex reflection group and $\mathfrak{h}$
its reflection representation over $\C$.  Let $\mathcal{S}$ denote
the set of complex reflections in $W$. Let $\omega$ be the canonical
symplectic form on $V = \h \oplus \h^*$. For $s\in \mathcal{S}$, let
$\omega_s$ be the skew-symmetric form that coincides with $\omega$
on $\im (\id_V - s)$ and has $\ker (\id_V - s)$ as its radical. Let
${\bf c}: \mathcal{S} \longrightarrow \C$ be a $W$-invariant
function sending $s$ to $c_s$. The {\it rational Cherednik algebra} at
parameter $t =0$ (depending on ${\bf c}$) is the quotient of the
skew group algebra of the tensor algebra on $V$, $TV \ast W$, by the
relations \begin{equation} \label{relation} [x,y] = \sum_{s\in
\mathcal{S}} c_s \omega_s(x,y) s\end{equation} for all $x,y\in V$.
This algebra is denoted by $H_{\bf c}$.

Let $Z_{\bf c}$ denote the centre of $H_{\bf c}$ and set $A =
\C[\h^*]^W \otimes \C[\h]^W$. Thanks to \cite[Proposition 4.15]{EG}
$A \subset Z_{\bf c}$ for any parameter ${\bf c}$ and $Z_{\bf c}$ is
a free $A$-module of rank $|W|$. Let $X_{\bf c}$ denote the spectrum
of $Z_{\bf c}$: this is called the {\it Calogero-Moser space} associated to $W$. Corresponding to the inclusion $A\subset Z_{\bf c}$
there is a finite surjective morphism $ \label{upsilon}
\Upsilon_{\bf c} : X_{\bf c} \longrightarrow \h^*/W \times \h/W.$

Let $\mathfrak{m}$ be the homogeneous maximal ideal of $A$. The {\it
restricted rational Cherednik algebra} is $\overline{H}_{\bf c} =
H_{\bf c}/\mathfrak{m}H_{\bf c}$. By \cite[PBW theorem 1.3]{EG} it
has dimension $|W|^3$ over $\C$. General theory asserts that the
blocks of $\overline{H}_{\bf c}$ are labelled by the closed points
of the scheme-theoretic fibre $\Upsilon^{*}(0)$. We call these points the {\it nilpotent points} of $X_{\bf c}$. By \cite[5.4]{G03}
there is a surjective mapping $$\Theta_{\bf c}: \irr W
\longrightarrow \{ \mbox{closed points of }\Upsilon_{\bf c}^{*}(0)
\} = \{ \mbox{blocks of }\overline{H}_{\bf c} \},$$ constructed by
associating to any $\blambda \in \irr W$ an indecomposable
$\overline{H}_{\bf c}$-module, the baby Verma module $M_{\bf
c}(\blambda)$. The fibres of $\Theta_{\bf c}$ partition $\irr W$. We will
call this the {\it $CM_{\bf c}$-partition}.

\subsection{} Let $W$ be a Weyl group. Let ${\bf L}: W \to \Q$ be a weight
function (in the sense of \cite[Section 2]{BGIL06}). Let
$\mathcal{H}$ be the corresponding Iwahori-Hecke algebra at unequal
parameters, an algebra over the group
algebra of $\Q$, $A = \oplus_{q\in \Q} \Z v^q$, which has a basis $T_w$
for $w\in W$, with multiplication given by the rule \[ T_sT_w =
\begin{cases} T_{sw}\ & \mathrm{if}\ l(sw)=l(w)+1\\ T_{sw}+
(v^{{\bf L}(s)}-v^{-{\bf L}(s)})T_w\ & \mathrm{if}\ l(sw)=l(w)-1
\end{cases} \] where $s\in \mathcal{S}$ and $w\in W$.
There is an associated partition of $W$ into
two-sided cells, see \cite[Chapter 8]{L03}. We call these the 
$KL_{\bf L}$-cells.

 \begin{conjecture} \label{theconj} Let $W$ be a Weyl group and let ${\bf L}$ be the weight
 function generated by ${\bf L}(s)=c_s$ for each $s\in \mathcal{S}$.
 \begin{enumerate}
 \item There is a natural identification of the $CM_{\bf c}$-partition and the $KL_{\bf L}$-partition; this is induced by attaching a $KL_{\bf L}$-cell to an irreducible $W$-representation via the asymptotic algebra $J$, \cite[20.2]{L03}.
  \item Let $\mathcal{F}$ be a $KL_{\bf L}$-cell of $W$ and let $M_{\mathcal{F}}$ be the closed point of $\Upsilon_{\bf c}^{*}(0)$ corresponding to $\mathcal{F}$ by (1). Then $\dim_{\C} (\Upsilon_{\bf c}^*(0)_{M_{\mathcal{F}} })= |\mathcal{F}|.$
 \end{enumerate}
 \end{conjecture}

The existence of the asymptotic algebra mentioned in (1) is still a conjecture, depending on Lusztig's conjectures P1-P15 in \cite[Conjecture 14.2]{L03}.

 \subsection{} This conjecture generalises the known results about
 the blocks of $\overline{H}_{\bf c}$ and about the fibre $\Upsilon_{\bf c}^{*}(0)$.
 \begin{itemize}
 \item \cite[Corollary 5.8]{G03} If $X_{\bf c}$ is smooth then
 $\Theta_{\bf c}$ is bijective, making the $CM_{\bf c}$-partition
 trivial. If $S\in \irr W$ then $\dim_\C (\Upsilon_{\bf
c}^*(0)_{M_{\mathcal{S}}}) = \dim_\C (S)^2$.
 \item $\Theta_{\bf c}$ is not bijective when $W$ is a finite
 Coxeter group of type $D_{2n}$, $E$, $F$, $H$ or $I_2(m)$ ($m\geq 5$),
 \cite[Proposition 7.3]{G03}. In all of these cases there are non-trivial two-sided cells.
 \item Both the $a$-function and the $c$-function are constant
 across fibres of $\Theta_{\bf c}$, \cite[Lemma 5.3 and Proposition 9.2]{G07}. This should be a property of two-sided cells.
 \end{itemize}
An advantage of the Cherednik algeras is that the $CM_{\bf c}$-partition exists for any complex reflection group whereas, at the moment, a cell theory only exists for Coxeter groups.

 \subsection{} In Theorem \ref{thethm} we will give a
 combinatorial description of the $CM_{\bf c}$-partition
 when $W = G(\ell , 1, n) = \mu_{\ell} \wr \mathfrak{S}_n$, and then in Theorem \ref{theweethm} we will provide evidence for the conjecture by showing that the $CM_{\bf c}$-partition agrees with the conjectural
description of the $KL_{{\bf L}}$-partition for $W = G(2, 1, n)$, the Weyl group
of type $B_n$, given in \cite[Section 4.2]{BGIL06}.

  \section{Blocks for $W = G(\ell, 1, n)$}

\subsection{} Let $\ell$ and $n$ be positive integers. Let $\mu_{\ell}$
be the group of $\ell$-th roots of unity in $\C$ with generator
$\sigma$ and let $\mathfrak{S}_n$ be the symmetric group on $n$
letters. Let $W$ be the wreath product $G(\ell ,1, n)= \mu_{\ell}\wr
\mathfrak{S}_n = (\mu_{\ell})^n\rtimes \mathfrak{S}_n$ acting naturally on $\h = \C^n$.

\subsection{} Let $\mathcal{P}(n)$ denote the set of partitions of
$n$ and $\mathcal{P}(\ell ,n)$ the set of $\ell$-multipartitions of $n$.
The set $\irr W$ can be identified naturally with
$\mathcal{P}(\ell,n)$ so that the trivial representation corresponds
to $((n), \emptyset, \ldots , \emptyset)$, e.g. \cite[Theorem 4.4.3]{JK81}. Given an element ${\bf s}
\in \mathbb{Z}_0^{\ell} = \{ (s_1, \ldots , s_{\ell})\in \Z^{\ell}: s_1 + \cdots
+ s_{\ell} = 0 \}$ there is an associated $\ell$-core (a
partition from which no $\ell$-hooks can be removed). The inverse of
the process assigning to a partition its $\ell$-core and
$\ell$-quotient defines a bijection \begin{equation} \label{multi}
\mathbb{Z}_0^{\ell} \times \coprod_n \mathcal{P}(\ell ,n)
\longrightarrow \coprod_n \mathcal{P}(n), \qquad ({\bf s}, \blambda)
\mapsto \tau_{\bf s} (\blambda).\end{equation} A detailed discussion
of this can be found in \cite[Section 2.7]{JK81} or \cite[Section
6]{G07}.

\subsection{} The Young diagram of a partition $\lambda$ will always
be justified to the northwest (one of the authors is English); we
will label the box in the $p$th row and $q$th column of $\lambda$ by
$s_{pq}$.  With this convention the residue of $s_{pq}$ is defined
to be congruence class of $p-q$ modulo $\ell$.  Recall that $s_{pq}$
is said to be $j$-{\it removable} for some $0\leq j \leq \ell-1$ if it has
residue $j$ and if $\lambda \setminus \{s_{pq}\}$ is the Young
diagram of another partition, a {\it predecessor} of $\lambda$. We say
that $s_{pq}$ is $j$-{\it addable} to $\lambda\setminus \{ s_{pq}\}$.

Let $J\subseteq\{0, \ldots , \ell -1\}$. We define the {\it $J$-heart}
of $\lambda$ to be the sub-partition of $\lambda$ which is obtained
by removing as often as possible $j$-removable boxes with $j\in J$
from $\lambda$ and its predecessors. 
A subset of $\mathcal{P}(n)$ whose elements are the partitions with
a given $J$-heart is called a {\it $J$-class}.

\subsection{} We will use the ``stability parameters" of \cite{G07} ${\btheta}({\bf c}) =
(\theta_0, \ldots , \theta_{\ell -1}) $ defined by $\theta_ k  =
-\delta_{0k}c_{(i,j)} + \sum_{t=1}^{\ell-1} \eta^{tk} c_{\sigma^t}$
for $0\leq k \leq \ell -1$, $\eta$ a primitive $l$-th root of unity
and an arbitrary transposition $(i,j)\in \mathfrak{S}_n$: they
contain the same information as ${\bf c}$.
Following \cite[Theorem 4.1]{G07} we set $\Theta = \{ (\theta_0, \ldots , \theta_{\ell -1}) \in \mathbb{Q}^{\ell} \}$ and $\Theta_1 = \{ {\btheta}\in \Theta : \theta_0 + \cdots + \theta_{\ell -1} = 1 \}.$

Let $\tilde{\mathfrak{S}}_\ell$ denote the affine symmetric group with generators
$\{\sigma_0, \dots ,\sigma_{\ell-1}\}$. It acts naturally on
 $ \Theta$ by $\sigma_j \cdot (\theta_0 , \ldots, \theta_{\ell-1}) =
 (\theta_0 , \ldots , \theta_{j-1}+\theta_j, -\theta_j, \theta_j+\theta_{j+1},
  \ldots , \theta_{\ell -1}).$ This restricts to the affine reflection representation
 on $\Theta_1$: the {\it walls} of $\Theta_1$ are
the reflecting hyperplanes and the {\it alcoves} of $\Theta_1$ are
the connected components of (the real extension of)
$\Theta_1\setminus \{ \text{walls}\}$. Let $A_0$ be the alcove
containing the point $\ell^{-1} (1, \ldots , 1)$: its  closure is a
fundamental domain for the action of $\tilde{\mathfrak{S}}_{\ell}$
on $\Theta_1$. The stabiliser of a point ${\btheta} \in
\overline{A}_0$ is a standard parabolic subgroup of
$\tilde{\mathfrak{S}}_{\ell}$ generated by simple reflections $\{
\sigma_j : j\in J\}$ for some subset $J\subseteq \{0, \ldots , \ell
-1\}$. We call this subset the {\it type} of $\btheta$. The type of
an arbitrary point $\btheta \in \Theta_1$ is defined to be the type
of its conjugate in $\overline{A}_0$.

\subsection{} We have an isomorphism $\tilde{\mathfrak{S}}_{\ell}
\cong \mathbb{Z}_0^{\ell}\times \mathfrak{S}_{\ell}$ with
$\mathfrak{S}_{\ell} = \langle \sigma_1, \ldots , \sigma_{\ell -1}\rangle$ and
the elements of $\mathbb{Z}_0^{\ell}$ corresponding to translations
by elements of the dual root lattice $\mathbb{Z}R^{\vee}$. The
symmetric group $\mathfrak{S}_{\ell}$ acts on $\mathcal{P}(\ell ,n)$
by permuting the partitions comprising an $\ell$-multipartition.

\begin{thm}
\label{thethm}Assume that $\btheta({\bf c}) \in \Theta_1$, so that
${\btheta}({\bf c})$ had type $J$ and belongs to $({\bf s}, w)\cdot \overline{A}_0$
for some $({\bf s}, w)\in \tilde{\mathfrak{S}}_{\ell}$. Then $\blambda , \bmu \in \irr W = \mathcal{P}(\ell ,n)$ belong
to the same block of $\overline{H}_{\bf c}$ if and only if
$\tau_{\bf s}(w\cdot \blambda)$ and $\tau_{\bf s}(w\cdot \bmu)$
belong to the same $J$-class. In other words, the $CM_{\bf
c}$-partition is governed by the $J$-classes.
\end{thm}

\begin{proof} Rescaling gives an isomorphism between $\overline{H}_{\bf c}$
and $\overline{H}_{{\bf c}/2}$ so we can replace $\bf c$ by ${\bf
c}/2$. By \cite[5.4]{G03} we must show that the baby Verma
modules $M_{{\bf c}/2}(\blambda)$ and $M_{{\bf c}/2}(\bmu)$ give
rise to the same closed point of $\Upsilon_{{\bf c}/2}^{\ast}(0)$ if
and only if $\tau_{\bf s} (w\cdot \blambda)$ and $\tau_{\bf s}
(w\cdot \bmu)$ have the same $J$-class. But the closed points of
$\Upsilon_{{\bf c}/2}^{\ast}(0)$ correspond to the $\C^*$-fixed
points of $X_{{\bf c}/2}$ under the action induced from the grading on
$H_{{\bf c}/2}$ which assigns degree $1$, respectively $-1$ and $0$,
to non-zero elements of $\h$, respectively $\h^*$ and $W$. By
\cite[Theorem 3.10]{G07} these agree with the $\C^*$-fixed points on
the affine quiver variety $\mathcal{X}_{\btheta ({\bf c})}(n)$ and
hence, thanks to \cite[Equation (3)]{G07} to the $\C^*$-fixed
points on the Nakajima quiver variety $\mathcal{M}_{\btheta ({\bf
c})}(n)$. Now the result follows since the combinatorial description
of these fixed points in \cite[Proposition 8.3(i)]{G07} is exactly
the one in the statement of the theorem.
\end{proof}

\subsection{Remarks}\label{remarks}
(1) The assumption $\btheta({\bf c}) \in \Theta_1$ imposes two
restrictions. First it places a
rationality condition on the entries of ${\bf c}$; guided by
corresponding results for Hecke algebras, \cite[Theorem
1.1]{dipmat}, we hope that this is not really a serious restriction.  Second it forces $c_{(i,j)} \neq 0$; if $c_{(i,j)} \neq 0$ then we can rescale to produce an
isomorphism $H_{\bf c} \cong H_{\lambda \bf c}$ and hence ensure
$\theta_0 + \cdots + \theta_{\ell -1} =1$.

(2) A generic choice of $\btheta({\bf c})\in \Theta_1$ will have
type $J = \emptyset$. The corresonding $CM_{\bf c}$-partition will
then be trivial and thus $X_{\bf c}$ will be smooth, \cite[Corollary
1.14(i)]{EG}.

\section{The case $W = G(2,1,n)$}

\subsection{}\label{chariso} We now focus on the situation where $W = G(2,1,n)$,
the Weyl group of type $B_n$. Here there are two conjugacy classes
of reflections $s$ and $t$, containing $(i,j)$ and $\sigma$
respectively. We will always assume that ${\bf c} = (c_s, c_t)\in
\mathbb{Q}^2$ has the property that $c_s,c_t\neq 0$. Corresponding
to the two group homomorphisms $\epsilon_1, \epsilon_2: W \to \C^*$,
$\epsilon_k(i,j)= (-1)^k$ for all $(i,j)\in s$ and
$\epsilon_k(\sigma)=(-1)^{k+1}$, there exist algebra isomorphisms
$H_{(c_s,c_t)} \cong H_{(-c_s,c_t)}$ and $H_{(c_s,c_t)} \cong
H_{(c_s,-c_t)}$, \cite[5.4.1]{GGOR}. So, without loss of
generality, we may assume that ${\bf c}\in \Q_{> 0}^2$.



\subsection{} There is a conjectural description of the
two-sided cells in \cite[Section 4.2]{BGIL06} which we recall very
briefly; more details can be found in both [loc.cit] and \cite{P07}.

We assume ${\bf L}(s) =a , {\bf L}(t) = b$ with $a,b \in \Q_{>0}$ and set $d= b/a$.
If $d\notin \mathbb{Z}$ then the partition is conjectured to be
trivial, \cite[Conjecture A(c)]{BGIL06}. If $d=r+1\in \mathbb{Z}$ then let $\mathcal{P}_r(n)$ be the
set of partitions of size $\frac{1}{2}r(r+1) + 2n$ with 2-core
$(r,r-1, \ldots , 1)$.  A domino tableau $T$ on $\lambda \in
\mathcal{P}_r(n)$ is a filling of the Young diagram of $\lambda$
with $0$'s in the 2-core and then $n$ dominoes in the remaining
boxes, each labelled by a distinct integer between $1$ and $n$ which
are weakly increasing both vertically and horizontally.  There is a
process called {\it moving through an open cycle} which leads to an
equivalence relation on the set of domino tableaux. This in turn
leads to an equivalence relation on partitions in $\mathcal{P}_r(n)$
where $\lambda$ and $\mu$ are related if there is a sequence of
partitions $\lambda = \lambda_0, \lambda_1, \ldots, \lambda_{s-1},
\lambda_s =\mu$ such that for each $1\leq i \leq s$, $\lambda_{i-1}$
and $\lambda_{i}$ are the underlying shapes of some domino tableaux
related by moving through an open cycle. The equivalence classes of
this relation are called $r$-cells. \cite[Conjecture D]{BGIL06}
conjectures that the two-sided cells are in natural bijection with
the $r$-cells.

\subsection{} The result of this section is the following.

\begin{thm} \label{theweethm} Under the bijection \eqref{multi} the $CM_{\bf c}$-partition
of $\irr W$ is identified with the above conjectural description of
the $KL_{{\bf L}}$-partition for ${\bf L}(s) = c_s, {\bf L}(t) = c_t$.
\end{thm}

This theorem shows that core-quotient algorithm provides a natural identification of the $CM_{\bf c}$-partition and the conjectural $KL_{\bf L}$-partition. We do not know in general whether Lusztig's conjectured mapping from $\irr W$ to the $KL_{\bf L}$-cells is given by this algorithm.

There are special cases where \cite[Conjecture
D]{BGIL06} has been checked -- for instance the asymptotic case $c_t
> (n-1)c_s$, \cite[Remark 1.3]{BGIL06} -- and thus in those cases
we really do get a natural identification between the $CM_{\bf c}$-classes and $KL_{\bf c}$-cells.

\subsection{} We will need the following technical lemma to prove the theorem.
\begin{lem} \label{thelem}
Let $\lambda\in \mathcal{P}_r(n)$ and set $j=r$ modulo $2$ with
$j\in \{0,1\}$. Suppose that $s_{pq}$ is a $j$-removable box and
$s_{tu}$ is a $j$-addable box such that $p\geq t$ and $q\leq u$ and
there are no other $j$-addable or $j$-removable boxes, $s_{vw}$,
with $p\geq v \geq t$ and $q\leq w \leq u$. Then there is a domino
tableau $T$ of shape $\lambda$ and an open cycle $c$ of $T$ such
that the shape of the domino obtained by moving through $c$ is
obtained by replacing $s_{pq}$ with $s_{tu}$.
\end{lem}

\begin{proof}
We use the notation of \cite[Sections 2.1 and 2.3]{P07} freely. We
consider the rim ribbon which begins at $s_{pq}$ and ends at
$s_{t,u-1}$. We claim that this rim ribbon can be paved by dominoes.
In fact this is a general property of a ribbon connecting a box,
$s$, of residue $j$ and a box, $e$, of residue $j+1$. Let $R$ be such
a ribbon. If $R$ contains only two boxes then $R = \{s,e\}$ so it is
clear. In general the box adjacent to $s$, say $s_{\mathrm{ad}}$,
has residue $j+1$ so that $R \setminus \{s, s_{\mathrm{ad}}\}$ is a
ribbon of smaller length and so the result follows by induction. In
our situation we can specify more. Starting at $s_{pq}$ we tile our
rim ribbon, $R$, as far as possible with vertical dominoes up to and
including $D= \{ s_{p-m+1,q}, s_{p-m,q}\}$ where $s_{p-m,q}$ has
residue $j+1$. If $s_{p-m,q}=s_{t,u-1}$ we have finished our tiling.
Otherwise $s_{p-m, q+1}\in R$ so the square $s_{p-m, q+1}$ will be
$j$-removable unless $\{ s_{p-m,q+1}, s_{p-m,q+2}\} \subseteq R$. We
now tile with as many horizontal dominoes as possible until we get to
$E= \{s_{p-m,q+k-1},s_{p-m,q+k}\}$ with $s_{p-m,q+k}$ having residue
$j+1$. If $s_{p-m,q+k} = s_{t,u-1}$ then our tiling stops. Otherwise
we must have the next domino as $F =
\{s_{p-m-1,q+k},s_{p-m-2,q+k}\}$ to avoid $s_{p-m,q+k+1}$ being
$j$-addable. We can now repeat this process to obtain our tiling of
$R$. From this description we obtain the following consequences. Let
$s_{vw}\in R$ have residue $j+1$ and suppose that $s_{vw}\neq
s_{t,u-1}$. Then
\begin{itemize}
\item[(i)] The domino which
contains $s_{vw}$ is either of the form $\{s_{v+1,w},s_{vw}\}$ or
$\{s_{v,w-1},s_{vw}\}$;
\item[(ii)] If $s_{v-1,w+1}\in \lambda$ then $s_{v,w+1} \in R$.
Furthermore, if $s_{vw}$ is contained in a horizontal domino then
$s_{v-1,w+1}\notin R$;
\item[(iii)] If $s_{v-1,w+1} \notin \lambda$ then $s_{v-1,w}\in R$.
\end{itemize}

Let $R$ denote the rim ribbon above and suppose it can be tiled by
$t$ dominoes. Let $\mu$ be the shape $\lambda \setminus R$. In
particular $\mu$ contains $\frac{1}{2}r(r+1) + 2(n-t)$ squares. By
the previous paragraph and \cite[Lemma 2.7.13]{JK81}, $\mu$ is a
Young diagram with the same $2$-core as $\lambda$ and so there
exists a $T^{\prime} \in \mathcal{P}_r(n-t)$ with shape $\mu$. Take
such a $T'$ filled with the numbers $1$ to $n-t$. Now add $R$ to
$T'$. We can tile $R$ by dominoes by the previous paragraph and we
fill the dominoes with the numbers $n-t+1, \ldots , n$ where the
filling is weakly increasing on the rows and columns of  $R$. This
gives a domino tableau $T = T' \cup R$ of shape $\lambda$.

We claim that $R \subseteq T$ is an open cycle, and that when we
move through this cycle we remove $s_{pq}$ from $T$ and add
$s_{tu}$. This will prove the lemma.

As we have seen in (i) a domino $D \subseteq R$ is either of the
form $D = \{ s_{vw}, s_{v+1,w} \}$ or $D=\{ s_{v,w-1}, s_{vw} \}$
with $s_{vw}$ having residue $j+1$. In the case that $D = \{ s_{vw},
s_{v+1,w} \}$ we have to study the square $s_{v-1, w+1}$ to
calculate $D'$. 
One of two things can happen. If this box does not belong to $T$
then $D' = \{s_{vw}, s_{v-1,w}\}$ and $D' \subseteq R$ by (iii).
Otherwise $s_{v-1, w+1}$ does belong to $T$. In this situation the
box is not in the rim so is filled with a lower value than $D$ and
so $D' = \{s_{vw}, s_{v,w+1}\}$. In particular, by (ii) above either
$D' \subseteq R$ or $D' =\{ s_{t,u-1},s_{tu}\}$.

Now suppose $D=\{ s_{v,w-1}, s_{vw} \}$. If the square $s_{v-1,
w+1}$ is not in $T$ then $D' = \{ s_{v-1,w}, s_{vw} \}$ and $D'
\subseteq R$ by (iii). If $s_{v-1, w+1}$ is in $T$ then it is not in
the rim by (ii) and so is filled with a value lower than that of
$D$. Thus $D'=\{ s_{vw}, s_{v,w+1} \}$ and $D'\subseteq R$ unless
$s_{vw}=s_{t,u-1}$, in which case $D'= \{ s_{t,u-1},s_{tu}\}$.

It is now clear $R$ that is a cycle and moving through this cycle
changes the shape of $\lambda$ by removing $s_{pq}$ and adding
$s_{tu}$.
\end{proof}
\subsection{Proof of Theorem \ref{theweethm}} We have ${\btheta}({\bf c}) = (-c_s + c_t,
-c_t)$ and by rescaling, see Remark \ref{remarks}(1), we consider
${\btheta}'({\bf c})=(1 - \frac{c_t}{c_s},\frac{c_t}{c_s}) \in
\Theta_1$. The action of $\tilde{\mathfrak{S}}_2$ on $\Theta_1$ is
given by $\sigma_0 \cdot (\theta_0,\theta_1) = (-\theta_0, \theta_1
+ 2\theta_0)$ and $\sigma_1 \cdot (\theta_0,\theta_1) = (\theta_0+2\theta_1,
-\theta_1)$. The walls are $\{(d,-d+1)\in \Theta_1: d\in
\Z\}$; they are of type $\{0\}$ if $d \in 2\Z$ and of type $\{1\}$
if $d\in 1+2\Z$. The fundamental alcove is $A_0 = \{ (d,-d+1): 0<d<1
\}$, and the alcove $A_r = \{ (d, -d+1) : r<d<r+1\}$ is then
labelled by either $((\frac{r}{2},\frac{-r}{2}), e)\in
\mathbb{Z}_0^2\times \mathfrak{S}_2$ or
$((\frac{-r+1}{2},\frac{r-1}{2}), \sigma_1)\in \mathbb{Z}_0^2\times
\mathfrak{S}_2$ depending on whether $r$ is even or odd.

If $\frac{c_t}{c_s}\notin \Z$ then ${\btheta}'({\bf c})$ has type
$\emptyset$ and the $CM_{\bf c}$-partition of $\irr W$ is trivial by
Theorem \ref{thethm} since $\emptyset$-classes are all singletons:
this agrees with the conjectured triviality of the two-sided cells
in this case. Thus we may assume
that $r = \frac{c_t}{c_s}-1 \in \Z_{\geq 0}$. Then ${\btheta}'({\bf
c})=(-r,r+1)$ will be in the closure of two alcoves, $A_{-r-1}$ and
$A_{-r}$. We consider the latter. Let ${\bf s}$ be the element in
$\Z_0^2$ coming from the labelling of $A_{-r}$; then$\tau_{\bf s}$ of \eqref{multi} produces a bijection between $\irr W
= \mathcal{P}(2,n)$ and $\mathcal{P}_r(n)$. If we set $j=r$
modulo $2$ and $J = \{ j\}$, the content of the theorem is simply
the assertion that the $r$-cells in $\mathcal{P}_r(n)$ consist of
the partitions in $\mathcal{P}_r(n)$ with the same $J$-heart.

Let us show first that if $\lambda,\mu\in\mathcal{P}_r(n)$ have the
same $J$-heart, say $\rho$, then they belong to the same $r$-cell.
The $J$-heart has no $j$-removable boxes, but we can construct the
partition $\mu$ from $\rho$ by adding, say, $t$ $j$-addable boxes.
Now let $\nu$ be the partition obtained from $\rho$ by adding $t$
$j$-addable boxes as far left as possible. We note that
by \cite[Theorem 2.7.41]{JK81} $\nu \in \mathcal{P}_r(n)$. Of course
$\mu$ and $\nu$ could be the same, but usually they will be
different. Now we apply Lemma \ref{thelem} again and again to $\nu$,
taking first the rightmost $j$-removable box from $\nu$ to the
position of the rightmost $j$-removable box on $\mu$ and then
repeating with the next $j$-removable box on the successor of $\nu$.
We continue until we have obtained a partition with shape $\mu$. By
Lemma \ref{thelem}, this process is obtained by moving through open
cycles. On the other hand, we can perform this operation  in the
opposite direction to move from $\lambda$ to $\nu$ via open cycles
(for this we use the same algorithm and the fact that for a cycle
$c$, moving through $c$ twice takes us back where we started,
\cite[Proof of Proposition 1.5.31]{G90}). It follows that $\lambda$
and $\mu$ belong to the same $r$-cell.

Finally, we need to see that if $\lambda,\mu \in \mathcal{P}_r(n)$
belong to the same $r$-cell, then they have the same $J$-heart. For
this it is enough to assume that $\mu$ is the shape of a tableau
obtained by moving through an open cycle on a tableau of shape
$\lambda$. But in this case the underlying shapes differ only in
some $j$-removable boxes, \cite[Section 2.3]{P07} and so they
necessarily have the same $J$-heart. \hfill $\Box$

\begin{ackn}
The second author gratefully acknowledges the support of The
Leverhulme Trust through a Study Abroad Studentship
(SAS/2005/0125).
\end{ackn}

\bibliographystyle{plain}
\bibliography{blocksRRCA}

\end{document}